\documentclass[12pt]{amsart}
\usepackage{latexsym}
\usepackage{amssymb}
\usepackage[dvips]{epsfig}
\usepackage{graphics,amsmath,amssymb}
\usepackage{amsmath}
\usepackage{amsfonts}
\usepackage[all]{xy}
\newtheorem{prop}{Proposition}[section]

\newtheorem{lemma}{Lemma}[section]

\title{Production matrices and Riordan arrays}

\author{Emeric Deutsch\and Luca
Ferrari\and Simone Rinaldi}

\address{Polytechnic University, Six
Metrotech Center, Brooklyn, NY 11201,\quad
{\tt deutsch@duke.poly.edu}\qquad (\textnormal{Emeric Deutsch}).}

\address{Dipartimento di Sistemi e Informatica, viale Morgagni 65, 50134 Firenze, Italy,
\quad {\tt ferrari@dsi.unifi.it}\qquad \textnormal{(Luca Ferrari)}.}

\address{Dipartimento di Scienze Matematiche e Informatiche ``R. Magari'', Pian dei Mantellini, 44,
53100 Siena, Italy,\quad {\tt rinaldi@unisi.it}\qquad \textnormal{(Simone Rinaldi).}}

\thanks{L. Ferrari and S. Rinaldi partially supported by
MIUR project: \emph{Linguaggi formali e automi: metodi, modelli e applicazioni }.}

\begin{document}

\maketitle

\begin{abstract}
We translate the concept of succession rule and the ECO method into matrix notation, introducing the concept of a
\emph{production matrix}. This allows us to combine our method with other enumeration techniques using matrices,
such as the method of Riordan matrices. Finally we treat the case of rational production matrices, i.e. those
leading to rational generating functions.
\end{abstract}

\medskip

\medskip

{\small {\bf Keywords}: ECO method, Production matrices, Riordan arrays.

\medskip

{\bf AMS 2000 Subject Classification}: 05A15, 05C38.}

\section{Introduction}

In~\cite{DFR} a new type of infinite matrices, the so called {\em production matrices}, was introduced and
studied. Such matrices are simply a new way to represent succession rules, which allows us to work with them
using algebraic methods.

\medskip

Succession rules were first introduced by Julian West in \cite{W1,W2}, and later were used as a formal tool for
the enumeration with the ECO method \cite{ECO} (ECO stands for Enumeration of Combinatorial Objects). In a word,
it is a constructive method to produce all the objects of a given class, according to the growth of a certain
parameter (the \emph{size}) of the objects. Basically, the idea is to perform a prescribed set of operations
(sometimes called {\em local expansion}) on each object of size $n$, thus constructing a set of objects of the
successive size. This construction should induce a partition of all the objects of any given size (that is, all
the objects of a given size are produced exactly once from the objects of immediately lower size, through the ECO
construction). If an ECO construction is sufficiently regular, then it is often possible to describe it using a
\emph{succession rule}. Intimately related to the concept of succession rule is the notion of \emph{generating
tree}, which is the most common way of representing a succession rule. The principal applications of ECO method
are enumeration \cite{enum,genfun}, random generation \cite{rangen}, or exhaustive generation \cite{exhgen,exgen}
of various combinatorial structures. For all these topics we refer the reader to the rich survey \cite{ECO}.


\bigskip

In \cite{DFR} we have proposed a translation of the concept of succession rule, and hence of the ECO method, into
matrix notation. This is achieved by introducing two different possibly infinite matrices strictly related to a
succession rule, namely its \emph{production matrix}, and its \emph{ECO matrix}. The main goal of our approach is
to provide a representation of succession rules more suitable for computations; indeed we define some operations
on production matrices in order to reproduce well known operations on the numerical sequences they represent.
This leads to the determination of the generating functions of such sequences, often more easily than it was
previously done by other methods (see \cite{ECO,genfun,ruop}).

%

In this paper we deepen the study of production matrices and ECO matrices, in particular trying to establish
relations with some kinds of infinite matrices that have been recently introduced and studied, known as {\em
Riordan arrays} \cite{altchar,R,Sh,SGWW,s}.

In Section~\ref{rpm} a comparison with the Riordan matrices method is established, thus completing an
investigation started in \cite{mv}. In Section~\ref{erpm} we outline the main results concerning {\em exponential
Riordan matrices} \cite{erm}, and then study this concept from the point of view of production matrices. In the
last section we consider finite production matrices, representing finite succession rules. Throughout the whole
paper, a huge amount of examples are described or simply sketched.

\medskip

Before ending this introduction, just a few words of explanation concerning some of the notations used in the
work. We have selected a more or less standard symbol for the generating functions of frequently occurring
sequences. Namely,

\medskip

\begin{itemize}
    \item $C=\frac{1-\sqrt{1-4z}}{2z}=1+z+2z^2 +5z^3 +14z^4 +42z^5 +132z^6 +429z^7 +\ldots$ for the Catalan numbers
(A000108);
    \item $M=\frac{1-z-\sqrt{1-2z-3z^2 }}{2z^2}=1+z+2z^2 +4z^3 +9z^4 +21z^5
+51z^6 +127z^7 +\ldots$ for the Motzkin numbers (A001006);
    \item $R=\frac{1-z-\sqrt{1-6z+z^2 }}{2z}=1+2z+6z^2 +22z^3 +90z^4 +394z^5
+1806z^6 +8558z^7 +\ldots$ for the large Schr\"oder numbers (A006318);
    \item $S=\frac{1+z-\sqrt{1-6z+z^2 }}{4z}=1+z+3z^2 +11z^3 +45z^4 +197z^5
+903z^6 +4279z^7 +\ldots$ for the small Schr\"oder numbers (A001003);
    \item $T=1+zT^3 =1+z+3z^2 +12z^3 +55z^4 +273z^5
1428z^6 +7752z^7 +\ldots$ for the ternary numbers (A001764);
    \item $F=\frac{1-\sqrt{1-4z}}{z(3-\sqrt{1-4z})}=1+z^2 +2z^3 +6z^4 +18z^5
+57z^6 +186z^7 +\ldots$ for the Fine numbers (A000957);
    \item $B= \frac{1}{\sqrt{1-4z}}=1+2z+6z^2 +20z^3 +70z^4 +252z^5
+924z^6 +3432z^7 +\ldots$ for the central binomial coefficients (A000984).
\end{itemize}

\medskip

Throughout the paper the A****** number between parentheses following a sequence is the identification number of
that sequence in \cite{EIS}. Most of the matrices we are going to consider are infinite; their lines (rows and
columns) will be indexed by nonnegative integers, and we will write ``line 0" to mean the first line, ``line 1"
to mean the second line, and so on.

\medskip

For the sake of readability we have decided that the paper be self-contained. This implies, in particular, that
some basic definitions and properties of production matrices, reported in Section 2, are taken directly from
\cite{DFR}. Therefore the reader familiar with these topics can skip Section 2.


\section{Basic definitions}\label{basics}

A succession rule is a formal system consisting of an \emph{axiom} $(a)$, $a\in \mathbf{N}^{+}$, and a set of
\emph{productions}:
\begin{displaymath}
\{ (k_{t})\rightsquigarrow (e_{1}(k_{t}))(e_{2}(k_{t}))\cdots
(e_{k_{t}}(k_{t})): t\in \mathbf{N}\} ,
\end{displaymath}
where $e_{i}:\mathbf{N}^{+}\longrightarrow \mathbf{N}^{+}$, which
explains how to derive the \emph{successors} $(e_{1}(k)),
(e_{2}(k)), \ldots , (e_{k}(k))$ of any given label $(k)$, $k\in
\mathbf{N}^{+}$. In general, for a succession rule $\Omega$, we
use the more compact notation:
\begin{equation}\label{regola}
\Omega :\left\{ \begin{array}{ll}
(a)
\\ (k)\rightsquigarrow (e_{1}(k))(e_{2}(k))\cdots (e_{k}(k))
\end{array}\right. .
\end{equation}

$(a)$, $(k)$, $(e_i(k))$, are called the {\em labels} of $\Omega$ (where $a,k,e_i (k)$ are positive integers).
The rule $\Omega$ can be represented by means of a {\em generating tree}, that is a rooted tree whose vertices
are the labels of $\Omega$; $(a)$ is the label of the  root and each node labeled $(k)$ has $k$ sons labeled by
$e_1(k), \ldots ,e_k(k)$ respectively, according to the production of $(k)$ in (\ref{regola}). A succession rule
$\Omega$ defines a sequence of positive integers $(f_n)_{n\geq 0}$, $f_n$ being the number of nodes at level $n$
in the generating tree determined by $\Omega$. By convention the root is at level $0$, so $f_0=1$. The function
$f_\Omega(x)=\sum _{n\geq 0} f_n x^n $ is the {\em generating function} determined by $\Omega$.

\bigskip

In \cite{DFR} we have proposed a representation of succession rules by means of infinite matrices \(P =
(p_{k,i})_{k,i\geq 0}\). Assume that the set of the labels of a succession rule is $\{ (l_k)\}_k$, and in
particular that $l_0$ is the label of the axiom. Then we define \(p_{k,i}\) to be the number of labels \(l_i\)
produced by label \(l_k\). We call \(P\) the {\em production matrix} of the succession rule. Observe that the
first row of a production matrix gives precisely the production of the axiom.

Labels do not occur explicitly in this matrix representation of succession rules. However, they are the row sums
of $P$. In particular, the label \(l_0\) of the axiom is the first row sum of $P$.

\bigskip

{\em Example.} To the succession rule
\begin{eqnarray}\label{central}
\left\{ \begin{array}{ll}
(2)
\\ (2k)\rightsquigarrow (2)^k(4)\cdots (2k)(2k+2)
\end{array}\right. ,
\end{eqnarray}
there corresponds the production matrix
\begin{equation}\label{matScr}
P = \begin{pmatrix} 1&1&0&0&0&0&\ldots
\\ 2&1&1&0&0&0&\ldots
\\ 3&1&1&1&0&0&\ldots
\\ 4&1&1&1&1&0&\ldots
\\ 5&1&1&1&1&1&\ldots
\\ \vdots&\vdots&\vdots&\vdots&\vdots&\vdots&\ddots
\end{pmatrix}.
\end{equation}

\bigskip

\begin{figure}[htb]
\centerline{\hbox{\psfig{figure=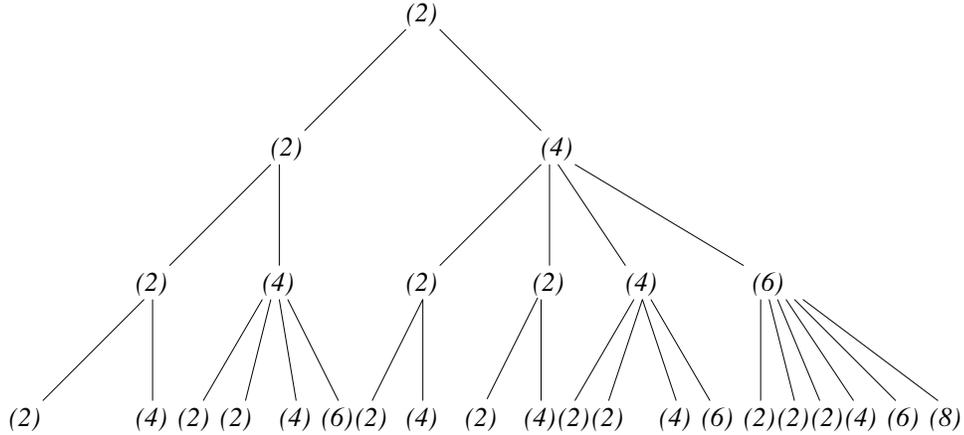,width=5in,clip=}}} \caption{The first levels of the generating tree
associated with the succession rule in~(\ref{central}).}
\end{figure}

\smallskip

In the generating tree at level zero we have only one node with
label \(l_0 (=2)\). This is represented by the row vector
\[
r_0 = \begin{pmatrix} 1&0&0&0&0&0&0&\ldots \end{pmatrix}.
\]

At the next levels of the generating tree the distribution of the labels  \(l_1,l_2,...\) is given by the row
vectors $r_i$, $i\geq 1$, defined by the recurrence relation $ r_i = r_{i-1}P$.

Stacking these row matrices, we obtain the matrix
\[
A_P=
\begin{pmatrix}
1&0&0&0&0&0 \quad \ldots \\ 1&1&0&0&0&0 \quad \ldots \\ 3&2&1&0&0&0 \quad \ldots \\
10&6&3&1&0&0
\quad \ldots \\ 35&20&10&4&1&0 \quad \ldots \\
126&70&35&15&5&1 \quad \ldots \\
\vdots &\vdots &\vdots &\vdots &\vdots &\vdots  \quad \ddots
\end{pmatrix}.
\]

The row sums of the above matrix are
\[
1,2,6,20,70,252,924,3432,12870,48620,\ldots
\]
i.e. the central binomial coefficients. This is the sequence corresponding to the succession rule of our example.
The enumerative properties of this succession rule have been examined in detail in \cite{BFR}. Using the
terminology of \cite{DFR}, we will refer to $A_P$ as the \emph{ECO matrix} induced by $P$.

\bigskip

\emph{Remarks.}\quad Let \(P \) be the production matrix of a given succession rule \(\Omega\). Throughout the
whole paper we will denote

$$ u^\top = (1,0,0,\ldots ,0,\ldots ) \qquad \qquad \qquad e= \begin{pmatrix}
1 \\
1 \\
1 \\
\vdots \\
\end{pmatrix}
$$

having appropriate sizes. The following properties have already been stated in~\cite{DFR}.
\begin{itemize}
\item[(i)] The labels of the nodes of the corresponding generating tree are the row sums of \(P\). If two row
sums happen to be equal, then, as labels, they will be considered to be distinct. This can be achieved by using,
for example, distinguishing subscripts; in the vocabulary of succession rules, these are called \emph{colored
succession rules} (see [FPPR]). \item[(ii)] \label{qq} The distribution of the nodes having various labels at the
various levels is given by the ECO matrix
\[
A_P =
\begin{pmatrix}
u^\top \\
u^\top P \\
u^\top P^2 \\
... \\
\end{pmatrix}
\]
(indeed, we have \(r_0 = u^\top, r_1 = r_0 P = u^\top P, r_2=r_1
P= u^\top P^2, \ldots \)).  The same fact can be expressed in a
concise way by the matrix equality
\begin{equation}\label{sti}
DA_P =A_P P,
\end{equation}
where $D=(\delta_{i,j+1})_{i,j\geq 0}$ ($\delta$ is the usual
Kronecker delta). In some sources \cite {tridiag,Sh} the matrix
$P$ is also called the {\em Stieltjes transform matrix of $A_P$}.
\item[(iii)] The sequence \(a_n\) induced by the succession rule
is given by \(a_n=u^\top P^n e\).
\item[(iv)] The bivariate generating function of the matrix
\(A_P\) is
\[
    G(t,z) = u^\top (I - zP)^{-1} \begin{pmatrix} 1 \\ t \\ t^2
    \\ t^3
\\ \vdots \\ \end{pmatrix}. \]
\item[(v)] The generating function of the sequence corresponding to the succession rule is
\[
    f_P(z) = u^\top (I - zP)^{-1} e.
\]


\item[(vi)] The exponential generating function of the sequence corresponding to the succession rule is
\[
    F_P(z) = u^\top \exp(zP) e .
\]
\end{itemize}

\bigskip

{\em Example.} We intend to find the sequence determined by the
production matrix
\begin{equation}\label{Bell}
P=
\begin{pmatrix}
0&1&0&0&0&\ldots \\
0&1&1&0&0&\ldots \\
0&0&2&1&0&\ldots \\
0&0&0&3&1&\ldots \\
0&0&0&0&4&\ldots \\
\vdots &\vdots &\vdots &\vdots &\vdots &\ddots
\end{pmatrix}.
\end{equation}

Denoting by \(P_n\) the upper left \(n\) by \(n\) submatrix of
\(P\), it is not difficult to compute the exponential of the
matrix $zP_n$, since it is an upper triangular matrix. The
eigenvalues of $P_n$ are easily seen to be the nonnegative
numbers $0,1,2,3,\ldots ,n-1$, each with multiplicity 1. Thus we have
immediately
\begin{displaymath}
\exp (zP_n)=C\exp (zD_n)C^{-1},
\end{displaymath}
where $D_n$ is the diagonalization of $P_n$, so that $\exp (zD_n)=(z\delta _{i,j}e^{iz})_{0\leq i,j\leq n}$, and
$C$ is a suitable invertible matrix. More precisely, a simple computation shows that $C$ is an upper triangular
matrix in which the $(i,j)$ entry has the form $\frac{1}{(j-i)!}$ (where obviously $i\leq j$). This implies that
also $C^{-1}$ is upper triangular and, for $i\leq j$, its $(i,j)$ entry has the form
$(-1)^{j-i}\frac{1}{(j-i)!}$. Now the computation of the first row of \(\exp(zP_n)\) is immediate, and we find
for it
\[
\begin{pmatrix}
1 & e^z-1 & \frac{1}{2!}(e^z-1)^2 & \frac{1}{3!}(e^z-1)^3 & \ldots
& \frac{1}{(n-1)!}(e^z-1)^{n-1} \\
\end{pmatrix}.
\]

Taking the sum of these entries and letting \(n \to \infty\),
for the exponential generating function induced by \(P\) we obtain
\[
    G_P(z) = e^{e^z-1}.
\]

The corresponding sequence is $1,1,2,5,15,52,203,876,\ldots$ (A000110; Bell numbers).

\section{Riordan production matrices}\label{rpm}

An infinite lower triangular matrix \(A\) is called a {\em Riordan
matrix} if its column \(k\) (\(k=0,1,2,...\)) has generating
function \(d(z)(zh(z))^k\), where \(d(z)\) and \(h(z)\) are formal
power series with \(d(0)\ne 0\). If, in addition, \(h(0)\ne 0\),
then \(A\) is said to be a {\em proper Riordan matrix}. We may
write \(A = (d(z),h(z))\).

Riordan matrices were first introduced in \cite{SGWW}, where it
has been proved that they constitute a natural way of describing
several combinatorial situations.

Proper Riordan matrices are characterized by a fundamental
property found by Rogers \cite{R} in 1978 and then examined
closely by Sprugnoli \cite{s}. In fact, for a given proper Riordan
matrix \(A = (d(z),h(z)) = (d_{n,k})_{n,k\geq 0}\), there exist unique
sequences \(\alpha = (\alpha_0, \alpha_1,\alpha_2,\ldots )\)
(\(\alpha_0 \ne 0\)) and  \(\zeta = (\zeta_0,
\zeta_1,\zeta_2,\ldots )\) such that
\begin{itemize}
\item[(i)] every element in column 0 can be expressed as a linear
combination of all the elements in the preceding row, the
coefficients being the elements of the sequence \(\zeta\), i.e.
\begin{equation}\label{col0}
    d_{n+1,0} = \zeta_0 d_{n,0} + \zeta_1 d_{n,1} +
    \zeta_2 d_{n,2} +\ldots
\quad ;
\end{equation}
\item[(ii)] every element \(d_{n+1,k+1}\), not lying in column 0
or row 0, can be expressed as a linear combination of the elements
of the preceding row, starting from the preceding column on, the
coefficients being the elements of the sequence \(\alpha\), i.e.
\begin{equation}\label{colk}
    d_{n+1,k+1} = \alpha_0 d_{n,k} + \alpha_1 d_{n,k+1} +
\alpha_2 d_{n,k+2} +\ldots \quad .
\end{equation}
\end{itemize}
Conversely, the existence of such sequences \(\alpha\) and
\(\zeta\) ensures that the matrix \(A\) is a proper Riordan
matrix.

The sequences \(\alpha\) and \(\zeta\) will be called the {\em \(\alpha\)-sequence} and the {\em
\(\zeta\)-sequence} of the Riordan matrix. By abuse of notation, by \(\alpha\), \(\zeta\) we shall denote also
the generating functions of these sequences.

The functions \(d(z),h(z),\zeta(z), \alpha(z)\) are connected by
the relations
\[
    h(z) = \alpha(zh(z)), \quad d(z) =
    \frac{d_{0,0}}{1-z\zeta(zh(z))}.
\]

It is also known \cite{s} that, under the assumption \(d(0)=h(0)\ne 0\) we
have \(d(z)=h(z)\) if and only if \(\alpha (z)=d_{0,0} +
z\zeta(z)\).

At this stage, it is natural to investigate the relationship between the theory of production matrices and the
theory of Riordan matrices. It turns out that, if the ECO matrix \(A_P\), induced by a production matrix \(P\),
is Riordan, then the matrix \(P\) has a very simple structure. First we complete a result just outlined in
Section \ref{basics}.

\begin{lemma} Let $P=(p_{i,k})_{i,k\geq 0}$ be an infinite
production matrix and let $A_P=(d_{n,k})_{n,k\geq 0}$ be its associated ECO matrix. Then, for any $k,n\geq 0$, we
have that $d_{0,k} = \delta _{0,k}$ and
\begin{equation}\label{recur}
d_{n+1,k}= d_{n,0}p_{0,k}+d_{n,1}p_{1,k}+d_{n,2}p_{2,k}+ d_{n,3}p_{3,k}+\ldots ,
\end{equation}
where $\delta$ is the usual Kronecker delta.

Conversely, if $A_P$ and $P$ are infinite matrices such that
relations (\ref{recur}) hold, then $A_P$ is the ECO matrix induced
by $P$.
\end{lemma}

\emph{Proof.} The first part of the theorem is essentially statement (ii) in the Remarks of Section 2. The
converse is easy.\qed

\begin{prop} Let \(P\) be an infinite production matrix and let
\(A_P\) be the matrix induced by \(P\). Then \(A_P \) is a Riordan
matrix if and only if \(P\) is of the form
\begin{equation}\label{prodrior}
P=
\begin{pmatrix}
\zeta_0 & \alpha_0&0&0&0 \quad \ldots \\ \zeta_1 & \alpha_1& \alpha_0&0&0 \quad \ldots \\ \zeta_2 &
\alpha_2 & \alpha_1 & \alpha_0 &0 \quad \ldots \\
\zeta_3 & \alpha_3 & \alpha_2 & \alpha_1 & \alpha_0 \quad \ldots \\ \zeta_4 & \alpha_4 & \alpha_3 & \alpha_2 &
\alpha_1 \quad \ldots \\ \zeta_5 & \alpha_5 & \alpha_4 & \alpha_3 & \alpha_2 \quad \ldots
\\ \vdots &\vdots &\vdots &\vdots &\vdots \quad \ddots
\end{pmatrix}.
\end{equation}
Moreover, columns 0 and 1 of the matrix \(P\) are the \(\zeta\)-
and \(\alpha\)-sequences, respectively, of the Riordan matrix
\(A_P\).
\end{prop}

{\em Proof.} If \(P\) is as in (\ref{prodrior}), then from the
relation \(r_n=r_{n-1}P\) between the rows of the matrix $A_P$
(see the above lemma) we can see at once that columns 0 and 1 of
\(P\) play the roles of the \(\zeta\)- and \(\alpha\)-sequences,
respectively. Consequently, by what has been said at the beginning
of this section, the matrix \(A_P\) is Riordan. Conversely, assume
that \(A_P\) is a Riordan matrix. Then equalities (\ref{col0}) and
(\ref{colk}) hold, and in view of the above lemma they translate
into the fact that the production matrix of $A_P$ must be as in
(\ref{prodrior}).\qed

\bigskip


The above proposition, formulated in terms of succession rules,
is the main result of \cite{mv}.

Because of the above property, a production matrix \(P\) having
the form (\ref{prodrior}) will be called a {\em Riordan production
matrix}.

In the case of a given Riordan production matrix \(P\), having \(
(\zeta_n)_{n\geq 0}\) and \( (\alpha_n)_{n\geq 0}\) as its first
two columns, one can easily determine the bivariate generating
function \(G(t,z)\) of the matrix \(A_P\) induced by \(P\) and
then, obviously, also the generating function \(f_P(z)\) of the
sequence induced by \(P\).

\begin{prop} Let \(P\) be a Riordan production matrix and let
\(\zeta(z)\) and \(\alpha(z)\) be the generating functions of its
first two columns, respectively. Then the bivariate generating
function \(G(t,z)\) of the matrix \(A_P\) induced by \(P\) and the
generating function \(f_P(z)\) of the sequence induced by \(P\)
are given by
\begin{equation}
    G_P(t,z) = \frac{d(z)}{1-tzh(z)},
    \quad f_P(z)=\frac{d(z)}{1-zh(z)},
\end{equation}
where \(h(z)\) is determined from the equation
\begin{equation}
    h(z)=\alpha(zh(z))
\end{equation}
and \(d(z)\) is given by
\begin{equation}\label{d}
    d(z)=\frac{1}{1-z\zeta(zh(z))}.
\end{equation}
\end{prop}

\emph{Proof.} The proof is straightforward from general facts of Riordan matrices theory. The only thing which is
worth remarking is that, in the case of ECO matrices, $d(0)=1$; this is the reason for which the numerator of the
fraction in the r. h. s. of (\ref{d}) is 1 (in general it should be $d(0)$).\qed

\bigskip

{\em Example.} Consider the Riordan production matrix
\[
P=
\begin{pmatrix}
3&1&0&0&0&0 \quad \ldots \\
7&3&1&0&0&0 \quad \ldots \\
15&7&3&1&0&0 \quad \ldots \\
31&15&7&3&1&0 \quad \ldots \\
63&31&15&7&3&1 \quad \ldots \\
127&63&31&15&7&3 \quad \ldots \\
\vdots &\vdots &\vdots &\vdots &\vdots &\vdots \quad \ddots
\end{pmatrix}.
\]
Note that the row sums of \(P\), i.e. the labels of the generating
tree, are the Eulerian numbers $4,11,26,57,120,\ldots$ (A000295). We
have \(\alpha(z)=\frac{1}{(1-z)(1-2z)}\) and \(\alpha-z\zeta =
1\). Then, recalling some of the results stated at the beginning
of this section, $A_P$ is the Riordan matrix $(d(z),h(z))$ such
that
\begin{equation}\label{ex1}
    h=\frac{1}{(1-zh)(1-2zh)}, \quad d=h.
\end{equation}
Let us denote by $f_P$ the generating function of the sequence
determined by $P$. From the last theorem we get
\begin{equation}\label{gf1}
 f_P=\frac{h}{1-zh}.
\end{equation}

Eliminating \(d\) and \(h\) from equations (\ref{ex1}),
(\ref{gf1}), we obtain \((1+zf_P)^3=f_P(1-zf_P)\). The
substitution \(K=z+z^2f_P\) leads to the equation
\(K^3=(K-z)(2z-K)\), which is the equation giving the generating
function for the number of noncrossing connected graphs. The
sequence corresponding to \(f_P\) starts $1,4,23,156,1162,\ldots$
(A007297).

\bigskip

{\em Example.} Consider the Riordan production matrix
\[
P= \begin{pmatrix}
0 & 1 & 0 & 0 & 0 & 0  & \ldots \\
2 & 1 & 1 & 0 & 0 & 0  & \ldots \\
11 & 3 & 1 & 1 & 0 & 0  & \ldots \\
54 & 12 & 3 & 1 & 1 & 0 & \ldots \\
272 & 55 & 12 & 3 & 1 & 1  & \ldots \\
1427 & 273 & 55 & 12 & 3 & 1 & \ldots \\
\vdots &\vdots &\vdots &\vdots &\vdots &\vdots &\ddots
\end{pmatrix},
\]
where in columns 1,2,3,... we have the ternary numbers
\(\frac{1}{2n+1}\binom{3n}{n}\) and in column 0 the same
numbers diminished by 1. Thus, \(\alpha(z) = T(z)\) and
\(\zeta(z)=\frac{T-1}{z}-\frac{1}{1-z}\), where the generating
function of ternary numbers \(T(z)\) satisfies \(T(z)=1+zT^3(z)\).
Consequently, \(h=\alpha(zh)=T(zh)=1+zhT^3(zh)=1+zh^4\), i.e.
\(h(z)\) is the generating function of the sequence
\(\frac{1}{3n+1}\binom{4n}{n}, n=0,1,2,...\). After some
elementary manipulations we obtain \(d=h(1-zh)\) and \(f_P(z)=
h(z)\). The sequence corresponding to \(f_P\) starts
1,1,4,22,140,969,$\ldots$ (A002293).

\bigskip

{\em Remark.} The last example can be easily generalized to
construct a production matrix inducing the sequence
\(\frac{1}{(p-1)n+1}\binom{pn}{n}, n=0,1,2,...\).

\bigskip

In the next table we show numerous other examples. The first two entries
in a row define a Riordan production matrix, the next two entries
show the bivariate generating function of the induced matrix
\(A_P\) and the generating function \(f_P\) of the induced
sequence, followed by the first few terms of the sequence and the
identification number in \cite{EIS}.

\bigskip

\hspace{-.4cm}\begin{tabular}{l|l|l|l|l|l} \\
\(\zeta(z)\) & \(\alpha(z)\) & \( G_P(t,z) \) & \(f_P(z)\) &
sequence & A-number \\ \hline
\(1\) & \(1\) &
\(\frac{1}{(1-z)(1-tz)}\) & \(\frac{1}{(1-z)^2}\) &
1,2,3,4,5,6,7,8,9,... & A000027 \\ \(\frac{z}{1-z^2}\) &
\(\frac{1}{1-z^2}\) & \(\frac{T(z^2)}{1-tzT(z^2)}\) &
\(\frac{T(z^2)}{1-zT(z^2)}\) & 1,1,2,3,7,12,30,55,... & A047749 \\
\(1\) & \(1+z\) & \(\frac{1}{1-(1+t)z}\) & \(\frac{1}{1-2z}\) &
1,2,4,8,16,32,64,... & A000079 \\ \(1\) & \(\frac{1}{1-z}\) &
\(\frac{1}{(1-z)(1-tzC)}\) & \(\frac{C}{1-z}\) &
1,2,4,9,23,65,197,... & A014137 \\ \(\frac{1}{1-z}\) & \(1+z\) &
\(\frac{(1-z)(1-2z)}{(1-3z+z^2)(1-z-tz)}\) &
\(\frac{1-z}{1-3z+z^2}\) & 1,2,5,13,34,89,233,... & A001519 \\
\(\frac{1}{(1-z)^2}\) & \(1\) & \(\frac{(1-z)^2}{(1-3z+z^2)(1-tz)}\) & \(\frac{1-z}{1-3z+z^2}\) &
1,2,5,13,34,89,233,... & A001519 \\ \(1+z\) & \(1+z+z^2\) & \(\frac{M}{1-tzM}\) & \(\frac{M}{1-zM} \) &
1,2,5,13,35,96,267,... & A005773 \\ \(\frac{z}{1-2z}\) & \(\frac{(1-z)^2}{1-2z}\) &
\(\frac{F}{1-tzF}\) & \(C \) & 1,1,2,5,14,42,132,... & A000108 \\
\(\frac{C}{1-z}\) & \(1\) & \(\frac{(1-z)C^2}{1-tz}\) &
\(\frac{C-1}{z}\) & 1,2,5,14,42,132,429,... & A000108 \\
\(\frac{1}{1-z}\) & \(\frac{1}{1-z}\) & \(\frac{C}{1-tzC}\) &
\(\frac{C-1}{z} \) & 1,2,5,14,42,132,429,... & A000108 \\
\(\frac{1-z-z^2}{(1-z)(1-2z)}\) & \(\frac{(1-z)^2}{1-2z}\) & \(\frac{(1-zF)C^2}{1-tzF}\) & \(\frac{C-1}{z} \) &
1,2,5,14,42,132,429,... & A000108 \\ \(\frac{C-1}{z}\) & \(1\) & \(\frac{1+B}{2(1-tz)}\) & \(\frac{1+B}{2(1-z)}
\) & 1,2,5,15,50,176,... & A024718 \\ \(1\) & \(\frac{1+z}{1-z}\) & \(\frac{1}{(1-z)(1-tzR)}\) & \(\frac{S}{1-z}
\) & 1,2,5,16,61,258,... & A104858 \\ \(1\) & \(\frac{1}{(1-z)^2}\) & \(\frac{1}{(1-z)(1-tzT^2)}\) &
\(\frac{T}{1-z} \) & 1,2,5,17,72,345,... & A104859 \\ \(\frac{1+z-z^2}{1-z}\) & \(\frac{1}{1-z}\) &
\(\frac{FC}{1-tzC}\) & \(\frac{F-1}{z^2} \) & 1,2,6,18,57,186,622,... & A000957\\ \(\frac{1}{(1-z)^2}\) &
\(\frac{1}{1-z}\) & \(\frac{B}{C(1-tzC)}\) & \(B \) & 1,2,6,20,70,252,924,... & A000984 \\ \(1+z\) & \((1+z)^2\)
& \(\frac{C}{1-tz(C-1))}\) & \(B \) & 1,2,6,20,70,252,924,... & A000984 \\ \(\frac{1}{1-2z}\) &
\(\frac{1-z}{1-2z}\) &
\(\frac{S}{1-tzS}\) & \(R \) & 1,2,6,22,90,394,... & A006318 \\
\(\frac{1}{1-z}\) & \(\frac{1+z}{1-z}\) & \(\frac{1+zR}{1-tzR}\) & \(R \) & 1,2,6,22,90,394,... & A006318 \\
\(0\) & \(\frac{2-z}{1-z}\) & \(\frac{1}{1-tz(1+R)}\) & \(R \) & 1,2,6,22,90,394,... & A006318 \\
\(\frac{1+z}{1-z}\) & \(\frac{1+z}{1-z}\) & \(\frac{1}{(1-zR)(1-tzR)}\) & \( S^2 \) & 1,2,7,28,121,550,... &
A010683\\ \(\frac{1}{(1-z)^2}\) & \(\frac{1}{(1-z)^2}\) & \(\frac{T}{1-tzT^2}\) & \(T^2 \) & 1,2,7,30,143,728,...
& A006013 \\ \(z\) & \(1+z+z^2\) & \(\frac{M}{(1+zM)(1-tzM)}\) & \(\frac{1}{\sqrt{1-2z-3z^2}} \) &
1,1,3,7,19,51,141,... & A002426 \\
\(2\) & \(\frac{1}{1-z}\) & \(\frac{1}{(1-2z)(1-tzC)}\) &
\(\frac{C}{1-2z} \) & 1,3,8,21,56,154,440,... & A014318 \\

\(\frac{2-2z+z^2}{1-z}\) & \(\frac{1}{1-z}\) &
\(\frac{C}{(1-z)(1-tzC)}\) & \(\frac{C^2}{1-z} \) &
1,3,8,22,64,196,625,... & A014138 \\

\(\frac{2}{1-z^2}\) & \(\frac{1}{1-z}\) & \(\frac{2C-1}{1-tzC}\) &
\(\frac{2(C-1)}{z}-C \) & 1,3,8,23,70,222,726,... & A000782 \\

\(\frac{z}{1-z}\) & \(\frac{1}{1-z}\) & \(\frac{F}{1-tzC}\) &
\(\frac{C-F}{z} \) & 1,3,8,24,75,243,808,... & A000958 \\

\(\frac{2-z}{1-z}\) & \(\frac{1}{1-z}\) & \(\frac{C^2}{1-tzC}\) &
\(C^3\) & 1,3,9,28,90,297,1001,... & A000245 \\

\(\frac{2}{1-z}\) & \(\frac{1}{1-z}\) & \(\frac{B}{1-tzC}\) & \(
BC \) & 1,3,10,35,126,462,... & A001700 \\

\(2+z\) & \((1+z)^2\) & \(\frac{C^2}{1-tzC^2}\) & \( BC \) &
1,3,10,35,126,462,... & A001700 \\

\end{tabular}

\hspace{-.8cm}\begin{tabular}{l|l|l|l|l|l} \\
\(\zeta(z)\) & \(\alpha(z)\) & \( G_P(t,z) \) & \(f_P(z)\) & sequence & A-number \\ \hline

\(\frac{2-3z-z^2}{(1-z)(1-2z)}\) & \(\frac{(1-z)^2}{1-2z}\) & \(\frac{BC(1-zF)}{1-tzF}\) & \( BC \) &
1,3,10,35,126,462,... &
A001700 \\

\(\frac{1+C}{1-z}\) & \(1\) & \(\frac{(1-z)BC}{1-tz}\) & \( BC \)
& 1,3,10,35,126,462,... & A001700 \\

\(2(1+z)\) & \((1+z)^2\) & \(\frac{B}{1-tzC^2}\) & \(
\frac{B^2}{C} \) & 1,3,11,42,163,638,... & A032443 \\

\(\frac{1}{1-z}\) & \(\frac{2-z}{1-z}\) & \(\frac{S}{1-tz(1+R)}\)
& \( \frac{S-1}{z} \) & 1,3,11,45,197,903,... & A001003 \\

\(\frac{2}{1-z}\) & \(\frac{1+z}{1-z}\) & \(\frac{R}{1-tzR}\) & \(
\frac{S-1}{z} \) & 1,3,11,45,197,903,... & A001003 \\

\(\frac{2-z}{(1-z)^2}\) & \(\frac{1}{(1-z)^2}\) &
\(\frac{T^2}{1-tzT^2}\) & \( \frac{T-1}{z} \) &
1,3,12,55,273,1428,... & A001764 \\

\(\frac{C-1}{z}-\frac{1}{1-z} \) & \(C\) &
\(\frac{T(1-zT)}{1-tzT}\) & \( T \) & 1,1,3,12,55,273,... &
A001764 \\

\(\frac{2}{1-2z} \) & \(\frac{1}{1-2z}\) &
\(\frac{C(2z)}{1-tzC(2z)}\) & \( \frac{C(2z)}{1-zC(2z)}\) &
1,3,13,67,381,2307,... & A064062 \\

\(\frac{1}{1-z} \) & \(\frac{2}{1-z}\) &
\(\frac{1}{(1-zC(2z))(1-2tzC(2z))}\) & \( \frac{C(2z)}{1-zC(2z)}\)
& 1,3,13,67,381,2307,... & A064062 \\

\(\frac{3-2z}{1-z} \) & \(\frac{1}{1-z}\) & \(\frac{BC}{1-tzC}\) &
\( BC^2\) & 1,4,15,56,210,792,... & A001791 \\

\(\frac{1}{1-3z} \) & \(\frac{3}{1-3z}\) &
\(\frac{C(3z)}{1-tzC(3z)}\) & \( \frac{C(3z)}{1-zC(3z)} \) &
1,4,25,190,1606,... & A064063 \\

\end{tabular}

\bigskip

\section{Exponential Riordan production matrices}\label{erpm}

In this section we outline the main results concerning exponential Riordan
matrices as they are exposited in \cite{erm}, and then study this concept
from the point of view of production matrices. For the case of tridiagonal
matrices, see also \cite{A1,A2,tridiag}.

\bigskip

Let $d(z),h(z)$ be two formal power series such that $d(0)\neq 0\neq h(0)$. An \emph{exponential Riordan}
(briefly, {\it eR}) \emph{matrix} is an infinite lower triangular array $A=(a_{n,k})_{n,k\geq 0}$ whose column
$k$ ($k=0,1,2,\ldots$) has exponential generating function $C_{k}(z)=\frac{1}{k!}d(z)(zh(z))^{k}$. We will use
the notation $A=[d(z),h(z)]$ to denote the {\it eR} matrix determined by $d(z)$ and $h(z)$. From the definition
it follows at once that the bivariate generating function of $A$, namely
\begin{equation}
G_A(t,z) = \sum_{n,k}a_{n,k}t^k\frac{z^n}{n!},
\end{equation}
is given by
\begin{equation}
G_A(t,z)=d(z)exp(tzh(z)).
\end{equation}

Observe that $G_A(t,z)$ has been defined to be ordinary with respect to the variable $t$ and exponential with
respect to the variable $z$.

In \cite{erm} it is shown that {\it eR} matrices form a group (with respect to the usual multiplication
operation), whose identity element is $[1,1]$.

Looking at the above definition, one could expect many similarities with the theory of (classical) Riordan
matrices. From the point of view of the present work, one of the main analogies with the ordinary case is the
possibility of expressing every entry of an {\it eR} matrix as a linear combination of the elements of the
preceding row. More precisely, we have the following result \cite{erm}.

\begin{prop}\label{analog} Let \(A = (a_{n,k})_{n,k\geq 0}=[d(z),h(z)]\) be
an exponential Riordan matrix and let
\begin{equation}\label{c,r}
    c(y)=c_0+c_1y+c_2y^2+\ldots ,\quad r(y)=r_0+r_1y+r_2y^2+\ldots
\end{equation}
be two formal power series such that
\begin{equation}\label{system}
\left\{ \begin{array}{ll}
r(zh(z)) = (zh(z))'
\\ c(zh(z)) = \frac{d'(z)}{d(z)}
\end{array} \right. .
\end{equation}

Then
\begin{equation}\label{(i)}
    (i) \quad a_{n+1,0} = \sum_i i! c_i a_{n,i}
\end{equation}
\begin{equation}\label{(ii)}
    (ii) \quad a_{n+1,k}=r_0a_{n,k-1} + \frac{1}{k!}\sum_{i\geq
k}i!(c_{i-k}+kr_{i-k+1})a_{n,i}
\end{equation}
or, defining \(c_{-1}=0\),
\begin{equation}\label{(iii)}
    a_{n+1,k}=\frac{1}{k!} \sum_{i \geq k-1} i!(c_{i-k} +
k r_{i-k+1})a_{n,i}.
\end{equation}

Conversely, starting from the sequences defined by (\ref{c,r}), the
infinite array $(a_{n,k})_{n,k\geq 0}$ defined by (\ref{(iii)}) is
an exponential Riordan matrix.
\end{prop}

\emph{Remark.} The sequences $(c_{n})_{n\geq 0},(r_{n})_{n\geq 0}$
are called respectively the \emph{c-sequence} and the
\emph{r-sequence} of $A$. There is a clear analogy with the
$\alpha$-sequence and the $\zeta$-sequence of the classical case.
However, whereas for an ordinary Riordan matrix the coefficients
of the linear combinations in (\ref{colk}) do not depend on the
column index, in the exponential case they do, as it is clear from
formula (\ref{(iii)}).

\bigskip

Consider the infinite matrix $P$ of the coefficients in (\ref{(iii)})
\begin{equation}
P = \begin{pmatrix}
c_0         & r_0               & 0         & 0     & 0     & \ldots \\
1!c_1       & \frac{1!}{1!}(c_0+r_1)    & r_0       & 0 & 0 & \ldots \\
2!c_2       & \frac{2!}{1!}(c_1+r_2)    & \frac{2!}{2!}(c_0+2r_1)   & r_0   & 0 & \ldots \\
3!c_3       & \frac{3!}{1!}(c_2+r_3)    & \frac{3!}{2!}(c_1+2r_2)   & \frac{3!}{3!}(c_0+3r_1)   & r_0 & \ldots \\
 4!c_4      & \frac{4!}{1!}(c_3+r_4)    & \frac{4!}{2!}(c_2+2r_3)   & \frac{4!}{3!}(c_1+3r_2)   & \frac{4!}{4!}(c_0+4r_1) & \ldots \\
\vdots & \vdots & \vdots & \vdots & \vdots & \ddots
\end{pmatrix}
\end{equation}
or, in a compact form,
\begin{equation}
    P = (p_{i,j})_{i,j\geq 0}, \qquad p_{i,j}=\frac{i!}{j!}(c_{i-j} + jr_{i-j+1}) \quad (c_{-1}=0),
\end{equation}
so that formula (\ref{(iii)}) can be rewritten as a
matrix equality:
\[
AP=DA,
\]
where, as in formula (\ref{sti}), $D=(\delta_{i,j+1})_{i,j\geq 0}$ ($\delta$ is the usual Kronecker delta). The
matrix $P$ is ``almost'' lower triangular, meaning that $p_{i+1,i}=r_{0}$; whereas $p_{i+k,i}=0$ for $k\geq 2$.

Denoting by diag$(n)$, $n\geq -1$, the sequence $(p_{n+k,k})_{k\geq 0}$, the
matrix $P$ is characterized by the following facts:
\begin{itemize}
\item diag(-1) is constant \(r_0, r_0, r_0,\ldots\),
\item diag(0) is an arithmetic progression with first term $c_0$ and ratio \(r_1\),
\item diag(1), after division by \(1,2,3,\ldots\), is an arithmetic progression with first term $c_1$ and ratio \(r_2\),
\item diag(2), after division by \(1\cdot 2,2\cdot 3,3\cdot 4,\ldots\), is an arithmetic progression with first term \(c_2\) and ratio \(r_3\), etc.
\end{itemize}

If the matrix $P$ has nonnegative integer entries, then it can be viewed as a
production matrix. In this case it will be called an \emph{exponential
Riordan production matrix}. Clearly, in order that $[d(z),h(z)]$ be an ECO
matrix, we have to set $d(0)=1$ as the remaining initial condition of the
system (\ref{system}).

\bigskip

It is convenient to have an expression for the bivariate exponential generating
function of the matrix $P$, namely
\begin{equation}\label{bgf}
\varphi_P (t,z)=\sum_{n,k}p_{n,k}t^k \frac{z^n}{n!}.
\end{equation}

Observe that again $\varphi_P$ has been defined to be ordinary with respect
to the variable $t$ and exponential with respect to the
variable $z$. By simply replacing the values of
the $p_{n,k}$'s in (\ref{bgf}) we immediately obtain:
\begin{eqnarray*}
\varphi_P (t,z)&=&\sum_{n,k}\frac{n!}{k!}(c_{n-k}+kr_{n-k+1})t^k \frac{z^n}{n!}
\\ &=&\sum_{n,k}\frac{t^k z^k}{k!}c_{n-k}z^{n-k}+\sum_{n,k}\frac{t^k
z^k}{k!}kr_{n-k+1}z^{n-k}
\\ &=&\sum_{k}\frac{t^k z^k}{k!}\cdot \sum_{n}c_n z^n +t\sum_{k}\frac{t^k
z^k}{k!}\cdot \sum_{n}r_n z^n
\end{eqnarray*}
i.e.
\begin{equation}\label{Pbegf}
\varphi_P (t,z)=e^{tz}(c(z)+tr(z)).
\end{equation}

Clearly, setting $t=1$ gives the exponential generating function of the row sums of $P$. If $P$ is an {\it eR}
production matrix, these are the \emph{labels} of the succession rule induced by $P$, and their exponential
generating function is:
\[
\varphi_P (1,z)=e^z (c(z)+r(z)).
\]

We consider in detail two examples. In the first one we are given an {\it eR} matrix $A$ and we will look for a
production matrix that induces $A$, while in the second example we start with an exponential Riordan production
matrix $P$ and we will find the induced ECO matrix.

\bigskip

{\em Example.} Consider the infinite matrix of the unsigned Stirling numbers
of the first kind
\[
A=
\begin{pmatrix}
1 & 0 & 0 & 0 & 0 & 0  & \ldots \\
1 & 1 & 0 & 0 & 0 & 0  & \ldots \\
2 & 3 & 1 & 0 & 0 & 0  & \ldots \\
6 & 11 & 6 & 1 & 0 & 0 & \ldots \\
24 & 50 & 35 & 10 & 1  & 0 & \ldots \\
120 & 274 & 225 & 85 & 15 & 1 & \ldots \\
\vdots & \vdots & \vdots & \vdots & \vdots & \vdots & \ddots
\end{pmatrix}.
\]

Its bivariate exponential generating function is
$\frac{1}{1-z}\exp(-t \log(1-z))$ and, consequently, $d(z)=\frac{1}{1-z}$ and
$h(z)=-\frac{\log(1-z)}{z}$. Now from the
system (\ref{system}) we obtain $r(-\log(1-z))=c(-\log(1-z))=\frac{1}{1-z}$,
leading to $r(y)=c(y)=\exp(y)$. In order to identify the corresponding matrix
$P$, we find from (\ref{Pbegf}) its bivariate exponential generating function
$\varphi_P(t,z)=(1+t)\exp((1+t)z)$. Now it follows at once that $P$ is the
matrix of the binomial coefficients:
\[
P=
\begin{pmatrix}
1 & 1 & 0 & 0 & 0 & 0 & \ldots \\
1 & 2 & 1 & 0 & 0 & 0 & \ldots \\
1 & 3 & 3 & 1 & 0 & 0 & \ldots \\
1 & 4 & 6 & 4 & 1 & 0 & \ldots \\
1 & 5 & 10 & 10 & 5 & 1 & \ldots \\
1 & 6 & 15 & 20 & 15 & 6 & \ldots \\
\vdots & \vdots & \vdots & \vdots & \vdots & \vdots & \ddots
\end{pmatrix}
\]
and the labels are $2,4,8,16,32,\ldots$.

\bigskip

\emph{Example.} Consider the matrix of the falling factorials (permutation coefficients):
\[ P = \begin{pmatrix}
1 & 1 & 0 & 0 & 0 & 0 & \ldots
\\ 2 & 2 & 1 & 0 & 0 & 0 & \ldots
\\ 6 & 6 & 3 & 1 & 0 & 0 & \ldots
\\ 24 & 24 & 12 & 4 & 1 & 0 & \ldots
\\ 120 & 120 & 60 & 20 & 5 & 1 & \ldots
\\ \vdots & \vdots & \vdots & \vdots & \vdots & \vdots & \ddots
\end{pmatrix}.
\]

The labels (i.e. the row sums of $P$) represent the total number of
arrangements of a set with $n$ elements (A000522).

Looking at the diagonals of $P$, it follows at once that $P$ is an {\it eR} production matrix; namely, the one
determined by
\[
c(y)=\frac{1}{(1-y)^2}\quad ,\qquad r(y)=\frac{1}{1-y}.
\]

System (\ref{system}) becomes
\[
(zh(z))'=\frac{1}{1-zh(z)}, \qquad \frac{d'(z)}{d(z)} = \frac{1}{(1-zh(z))^2},
\]
with initial conditions $h(0)=1$, $d(0)=1$ and we obtain

\[
d(z)=\frac{1}{\sqrt{1-2z}}\quad, \qquad h(z)=\frac{1-\sqrt{1-2z}}{z},
\]
so that $A_P =\left[ \frac{1}{\sqrt{1-2z}},\frac{1-\sqrt{1-2z}}{z}\right]$. Using the definition of an {\it eR}
matrix, $A_P$ turns out to be
\[
A_P =\begin{pmatrix}
1 & 0 & 0 & 0 & 0 & 0 & \ldots
\\ 1 & 1 & 0 & 0 & 0 & 0 & \ldots
\\ 3 & 3 & 1 & 0 & 0 & 0 & \ldots
\\ 15 & 15 & 6 & 1 & 0 & 0 & \ldots
\\ 105 & 105 & 45 & 10 & 1 & 0 & \ldots
\\ 945 & 945 & 420 & 105 & 15 & 1 & \ldots
\\ \vdots & \vdots & \vdots & \vdots & \vdots & \vdots & \ddots
\end{pmatrix}.
\]

The rows of $A_P$ are the coefficients of the Bessel polynomials with exponents in decreasing order (see A001497;
the sequence of the row sums of $A_P$ is A001515).

\bigskip

In closing this section we remark that it is easy to pass from an ordinary
Riordan matrix $(d(z),h(z))$ to the exponential Riordan matrix defined by the
same pair of formal power series, that is $[d(z),h(z)]$. Indeed, if
$(d(z),h(z))=(a_{n,k})_{n,k\geq 0}$ and $[d(z),h(z)]=(\alpha _{n,k})_{n,k\geq
0}$, then we have
\[
\alpha_{n,k}=\frac{n!}{k!}a_{n,k}.
\]

\section{Finite and rational production matrices}

A succession rule is said to be {\em finite} if it has a finite number of productions. This property simplifies
the structure of the rule, and many problems still open for the whole class of succession rules turn out to be
solvable for the class of finite succession rules. In \cite{BDPR} it is proved that the {\em equivalence
problem}, i.e. the problem of establishing if two given rules define the same numerical sequence, is decidable
for finite succession rules. Moreover finite succession rules lead to rational generating functions, and have a
combinatorial interpretation in terms of regular languages \cite{rational}. More specifically, defining

\medskip

\begin{enumerate}
    \item $\mathcal R$ the set of rational generating functions of integer sequences ($Z$-rational functions, in the
notation of \cite{SS});
    \item ${\mathcal R}^+$ the set of rational generating functions of positive integer sequences;
    \item $REG$ the set of generating functions of regular languages;
    \item $\mathcal S$ the set of rational generating functions of succession rules;
    \item $\mathcal F$ the set of generating functions of finite succession rules,
\end{enumerate}

\medskip

\noindent from \cite{SS} and \cite{oper} we obtain the following inclusions:

\bigskip

\centerline{
\xymatrix{ %
& {REG}\ar @{}[dr]|{\rotatebox{-45}{$\subset$}}\\
\mathcal{F} \ar@{}[ur]|-{\rotatebox{45}{$\subset$}} \ar@{}[dr]|-{\rotatebox{-45}{$\subseteq$}}& \ar@{}[r]|-{ }
&\mathcal{R}^+ \ar@{}[r]|-{\subset}
&\mathcal{R}\\
& \mathcal{S} \ar@{}[ur]|-{\rotatebox{45}{$\subset$}} %
}}

\bigskip

The classes $\mathcal R$, $REG$, and $\mathcal F$ are decidable, while ${\mathcal R}^+$ is not decidable. It is
also conjectured that ${\mathcal F}={\mathcal S}$, i.e. every rational succession rule (which is a succession
rule having a rational g.f.) is equivalent to a finite one.

\medskip

One of the most popular finite succession rules is the following, defining Fibonacci numbers, $1,1,2,3,5,8,13,21,
\ldots$ (sequence A000045):

\begin{equation}\label{fifib}
\left \{
\begin{array}{l}
(1) \\ (1) \leadsto (2) \\ (2) \leadsto (1)(2).
\end{array}
\right.
\end{equation}

A production matrix associated with a finite succession rule has finite dimensions: the number of rows is equal
to the number of productions of the rule, whereas the number of columns is equal to the number of labels. For
example, the production matrix associated with the rule in (\ref{fifib}) is

\[
\begin{pmatrix} 0&1\\1&1\\
\end{pmatrix}.
\]

We recall that this example was also considered in \cite{W1}.

\medskip

In this section we will determine some properties of finite and rational production matrices, re-obtaining and
sometimes extending some results in \cite{genfun}, \cite{BDPR}, and \cite{oper}.

\begin{prop}
If \(P\) is a finite production matrix, then the induced generating function \(f_P(z)\) is rational.
\end{prop}

{\em Proof.} This follows at once from \(f_P(z) = u^\top (I-zP)^{-1} e\).

\begin{prop} Let \(P\) be a finite production matrix and let
\[
    t^q + c_1t^{q-1} + ... + c_{q-1}t + c_q
\]
be the minimal polynomial of \(P\). Then the sequence \((a_n)\) induced by \(P\) satisfies the linear homogeneous recurrence equation
\[
    a_n + c_1a_{n-1} + ... + c_{q-1}a_{n-q+1} + c_q a_{n-q}=0.
\]
\end{prop}

{\em Proof.} The matrix \(P\) satisfies
\[
    P^q + c_1P^{q-1} + ... + c_{q-1}P + c_q I =0,
\]
from where, premultiplying by \(u^\top P^{n-q}\) and postmultiplying by \(e\), we obtain the conclusion of the theorem.

\medskip

{\em Remarks.} 1. Clearly, this last theorem implies again that the generating function induced by a finite production matrix is rational.

2. The last theorem can be slightly improved. Instead of the minimal polynomial one can take the \(
P\)-annihilator of \(e\), i.e. the unique monic polynomial which generates the ideal of all polynomials \(g(t)\)
such that \(g(P)e=0\). This is a divisor of the minimal polynomial of the matrix \(P\) (see, for example,
\cite{hk}, p. 228). However, even this need not be the recurrence equation of lowest order satisfied by the
sequence induced by \(P\).

{\em Example.} Consider the succession rule defined by the production matrix

\[
P=\begin{pmatrix} 2&1&1&0\\0&3&0&0\\0&1&2&1\\0&1&1&3\\
\end{pmatrix}.
\]

The minimal polynomial of \(P\) is \(t^4-10t^3+36t^2-55t+30\), yielding the recurrence equation of order four:

\[
    a_n - 10a_{n-1}+36a_{n-2}-55a_{n-3}+30a_{n-4}=0.
\]

However, the \(P\)-annihilator of \(e\) is \(t^3-8t^2+20t-15\), yielding the recurrence equation of order three:

\[
    a_n - 8a_{n-1}+20a_{n-2}-15a_{n-3} = 0.
\]

This equation is satisfied by the sequence of the \(i\)-th row sums of the matrices \(I,P,P^2,...\), where
\(i=1,2,3,4\). But, in this case, the sequence of the first row sums of the matrices \(I,P,P^2,...\), i. e. the
sequence induced by \(P\), satisfies the recurrence equation of order two:

\[
    a_n-5a_{n-1}+5a_{n-2}=0.
\]

(the polynomial \(t^2-5t+5\) is a divisor of \(t^3-8t^2+20t-15\)).

%
%
%

\medskip

In the last part of this section we will focus on rational production matrices, i.e. production matrices whose
induced generating function is rational. The conjectured equivalence between finite and rational succession rules
is the motivation for the present investigation. We will just state a simple result on rational production
matrices along with some examples. The rational production matrices considered in such examples are not finite;
in some cases the determination of an equivalent finite production matrix is a challenging problem.

\begin{prop} Let \(P\) be a production matrix. If the vectors
\(e, Pe, P^2e,P^3e,...\) are linearly dependent, then the induced generating function \(f_P(z)\) is rational.
\end{prop}

\emph{Proof.} The hypothesis says that

             $$\alpha_0 e + \alpha_1 P e + ... + \alpha_k P^k e = 0$$

for some $k \in \mathbf{N}$ and $\alpha_0 ,..., \alpha_k \in \mathbf{R}$. Premultiplying by $u^{\top} P^n$ we get

             $$\alpha_0 a_n + \alpha_1 a_{n+1} + ... + \alpha_k a_{n+k} = 0$$

for any $n \in \mathbf{N}$. Therefore $( a_n )_{n \geq 0}$ satisfies a linear recurrence, which implies that its
generating function is rational. \qed

\smallskip

{\em Remark.} The above theorem generalizes Proposition 2 of [BBDFGG]. Indeed, the linear dependence \(P^2e =
\alpha Pe + \beta e\) is equivalent to the ``affine function" assumption of Proposition 2 of [BBDFGG].

\smallskip

{\em Example.} We consider the succession rule
\[
    [(2)\{(k) \rightsquigarrow (2)^{k-2}(2+(k \mod 2))(k+1)\}]
\]
(see [BBDFGG]). Its production matrix is
\[
P=\begin{pmatrix}
1 & 1 & 0 & 0 & 0 & 0  &\ldots \\
1 & 1 & 1 & 0 & 0 & 0  &\ldots \\
3 & 0 & 0 & 1 & 0 & 0  &\ldots \\
3 & 1 & 0 & 0 & 1 & 0  &\ldots \\
5 & 0 & 0 & 0 & 0 & 1  &\ldots \\
5 & 1 & 0 & 0 & 0 & 0  &\ldots \\
\vdots &\vdots &\vdots &\vdots &\vdots &\vdots &\ddots \\
\end{pmatrix}
\]
Introducing the vectors
\[
v=\begin{pmatrix} 0 & 1 & 0 & 1 & 0 & 1 & 0 & 1 & 0 & 1 \: ...\\ \end{pmatrix}^\top,
\]
\[
w=\begin{pmatrix} 1 & 2 & 3 & 4 & 5 & 6 & 7 & 8 & 9 & 10 \: ...\\ \end{pmatrix}^\top,
\]
we can easily derive
\[
    Pe = e + w, \quad Pv = e, \quad Pw = e+v+2w,
\]
from where
\[
    Pe = e+w, \quad P^2e=2e+v+3w, \quad P^3e=6e+3v+8w.
\]
Eliminating \(v\) and \(w\) from the last three equalities, we obtain
\[
    P^3e-3P^2e+Pe-e=0.
\]
Consequently, the generating function induced by the matrix \(P\) is rational. This generating function can be
found by the standard procedure from the recurrence equation that follows from the previous relation, namely
\(a_n-3a_{n-1}+a_{n-2}-a_{n-3}=0\) and the initial values \(a_0=1, a_1=2, a_2=5\) (the first components of the
vectors \(e, Pe, P^2e\), respectively). One finds
\[
    f_P(z) = \frac{1-z}{1-3z+z^2-z^3}.
\]

It is easy to observe that the rational production matrix $P$ is equivalent to the following finite production
matrix:

\[
\begin{pmatrix}
1 & 1 & 0 \\
1 & 1 & 1 \\
2 & 1 & 1 \\
\end{pmatrix}
\]

{\em Remark.} In an entirely similar manner we can derive that the succession rule
\[
    [(2)\{(k) \rightsquigarrow (2)^{k-2}(3-(k \mod 2))(k+1)\}]
\]
(see [BBDFGG]) induces the generating function \(\frac{1+z-2z^2}{1-z-6z^2+2z^3}\).

\bigskip

{\em Example.} We consider the succession rule defined by the production matrix
\[
P=\begin{pmatrix}
1 & 1 & 0 & 0 & 0 & 0  &\ldots\\
3 & 0 & 1 & 0 & 0 & 0  &\ldots\\
6 & 0 & 0 & 1 & 0 & 0  &\ldots\\
10 & 0 & 0 & 0 & 1 & 0 &\ldots\\
15 & 0 & 0 & 0 & 0 & 1 &\ldots\\
21 & 0 & 0 & 0 & 0 & 0 &\ldots\\
\vdots &\vdots &\vdots&\vdots&\vdots&\vdots &\ddots\\
\end{pmatrix}
\]
where in the first column the \(i\)-th entry is
\(\frac{i(i+1)}{2}\). Denoting by \([c_i]\) a column vector with
\(i\)-th entry \(c_i\), we have
\[
P[1]=\left[\frac{1}{2}i^2 + \frac{1}{2}i + 1\right], \quad
P[i]=\left[\frac{1}{2}i^2 + \frac{3}{2}i + 1\right], \quad
P[i^2]=\left[\frac{3}{2}i^2 + \frac{5}{2}i + 1\right].
\]
By linearity we obtain
\[
Pe=\frac{1}{2}[i^2] + \frac{1}{2}[i] + e, \quad
P^2e=\frac{3}{2}[i^2] + \frac{5}{2}[i] + 2e, \quad
P^3e=\frac{9}{2}[i^2] + \frac{17}{2}[i] + 6e.
\]
Eliminating \([i]\) and \([i^2]\) from the last three relations, we obtain
\[
    P^3e - 4P^2e + 3Pe - e =0,
\]
leading to the recurrence equation
\[
    a_n - 4a_{n-1}+3a_{n-2}-a_{n-3} = 0.
\]
This equation with the initial values \(a_0=1, a_1=2, a_2=6\),
obtained from the first entries of \(e, Pe,\) and \(P^2e\),
respectively, yield the generating function
\[
    f_P(z)=\frac{(1-z)^2}{1-4z+3z^2-z^3}.
\]
The first ten terms of the sequence are 1,2,6,19,60,189,595,1873,5896,18560... (sequence A052544 in \cite{EIS}).

\medskip

{\em Remark.} If in the previous example, in the first column of \(P\) we take, more generally,  the \(i\)-th
term to be \(\alpha i^2+\beta i+\gamma \), then, in the same manner, for the recurrence equation we obtain
\[
    a_n-(\alpha + \beta + \gamma + 3)a_{n-1} - (\alpha -\beta
    -2\gamma-3)a_{n-2} - (\gamma + 1)a_{n-3}=0.
\]

\medskip

{\em Example.} We consider the succession rule defined by the production matrix
\[
P=\begin{pmatrix}
\alpha     & 1 & 0 & 0 & 0 & 0  &\ldots\\
\alpha     & 1 & 1 & 0 & 0 & 0  &\ldots\\
\alpha + 1 & 0 & 0 & 1 & 0 & 0  &\ldots\\
\alpha + 1 & 1 & 0 & 0 & 1 & 0  &\ldots\\
\alpha + 2 & 0 & 0 & 0 & 0 & 1  &\ldots\\
\alpha + 2 & 1 & 0 & 0 & 0 & 0  &\ldots\\
\vdots&\vdots&\vdots&\vdots&\vdots&\vdots &\ddots\\
\end{pmatrix},
\]
where \(\alpha\) is a nonnegative integer. It is easy to see that
\[
    Pe = (\alpha + 1)e + v,
\]
where \( v=\begin{pmatrix} 0 & 1 & 1 & 2 & 2 & 3 & 3 & 4 & 4 & 5 \: ... \end{pmatrix} ^\top \). We also have
\[
Pv = e + v.
\]
Eliminating \(v \) from the last two relations, we obtain
\[
P^2e - (\alpha + 2)Pe + \alpha e = 0,
\]
leading to the recurrence equation
\[
    a_n - (\alpha + 2)a_{n-1} + \alpha a_{n-2}=0.
\]
Making use of the initial conditions \(a_0=1, a_1=\alpha + 1\), we
obtain the generating function
\[
    f_P(z) = \frac{1-z}{1-(\alpha + 2)z + \alpha z^2}.
\]
For \(\alpha = 0\) we obtain the sequence \(1,1,2,4,8,16,...\)
(A000079), for \(\alpha = 1\) we obtain the odd-indexed Fibonacci
numbers \(1,2,5,13,34,89,233,...\) (A001519), and for \(\alpha =
2,3,4\) we obtain sequences A007052, A018902, A018903,
respectively.

\end{document}